\documentclass[11pt]{amsart}
\usepackage{amssymb}
\usepackage{graphicx}
\usepackage{amsmath, color}
\usepackage{algpseudocode}
\usepackage{psfrag}
\newcommand{\PP}{\mathbb{P}}
\newcommand{\EE}{\mathbb{E}}

\begin{document}

\title{Cascade sizes in a branching process with Gamma distributed generations}

\author{James Burridge*}

\date{\today}

\maketitle

\noindent \begin{center}  *{\scriptsize Department of Mathematics,
University of Portsmouth,
Portsmouth, UK. \textsf{james.burridge@gmail.com} }\\
\end{center}

\vspace{.6cm}
\begin{abstract}
We derive an exact expression for the probability density function of the cascade size (total progeny) in a continuous state branching process when the generations are Gamma distributed. The distribution has application in the modelling of cascade processes such as landslides and electrical network failures.

\end{abstract}

\section{Continuous state branching processes}

\subsection{General Theory}

We may define a continuous state branching process in the following way, see \cite{Dob06,Har63} for further details. Without loss of generality we let the size of the zeroth generation be $X_0=1$. The size of the first generation is then drawn from some distribution $G$ with support $[0,\infty]$, so $X_1 \sim G$. We now write $G^{\ast n}$ for the distribution of the sum of $n$ independent copies of $X_1$. We extend this definition to non integral $n$ as follows; if $\tilde{g}$ is the Laplace transform of the density function for the distribution $G$, then $(\tilde{g})^n$ is the Laplace transform of the density function for $G^{\ast n}$. The size of the $n$th generation is then a random variable with the following distribution:
\begin{equation*}
X_n \sim G^{\ast X_{n-1}}.
\end{equation*}
This defines a branching process.

\subsection{A Gamma Branching Process}

The Gamma distribution $\Gamma(k,\theta)$ has the probability density function
\begin{equation*}
f(x) = \frac{x^{k-1}e^{-\frac{x}{\theta}}}{\Gamma(k) \theta^k}.
\end{equation*}
If $X \sim \Gamma(k,\theta)$ then $\EE(X)=k\theta$, and $Var(X)=k \theta^2$. The distribution has the following property which holds for all $n\in \mathbb{R}^+$:
\begin{equation*}
\Gamma(k,\theta)^{\ast n} \equiv \Gamma(nk,\theta).
\end{equation*}
We make use of this property in setting up the following  branching process. Let $X_0=1$ be the size of the zeroth generation, and let:
\begin{equation*}
X_1 \sim \Gamma(2,p).
\end{equation*}
We have made the choice $k=2$ here but the following analysis may be generalised to arbitrary $k$. The size of the $n$th generation is distributed as:
\begin{equation*}
X_n \sim \Gamma(2,p)^{\ast X_{n-1}} \equiv \Gamma(2 X_{n-1},p).
\end{equation*}
The total cascade size in an infinite system:
\begin{equation*}
Z= \sum_{k=0}^\infty X_k
\end{equation*}
may be infinite. In what follows we will compute the probability of this event: $\PP\{Z=\infty\}$.

\section{The cascade distribution}

It is possible to derive the exact probability density function for the total size, $Z$ of the cascade by taking the continuum limit of a discrete branching process.

\subsection{Negative binomial approximation to the Gamma distribution}

We begin by noting that the negative binomial distribution, which has probability mass function:
\begin{equation*}
b(n,r,q) = \frac{\Gamma(n+r)}{n! \Gamma(r)} (1-q)^r q^n
\end{equation*}
provides an arbitrarily close discrete approximation to the gamma distribution for appropriate choice of the parameters $r$ and $q$. The approximation is set up in the following way. We divide $[0,\infty]$ into a discrete lattice of constant spacing $\delta$, and let $X_{\delta}$ be a discrete random variable which approximates $X \sim \Gamma(k,\theta)$. Let $X_{\delta}$ have the probability mass function:
\begin{equation*}
\PP(X_\delta=n \delta) = b(n,r,q).
\end{equation*}
The mean and variance of $X_\delta$ are then:
\begin{align*}
\EE(X_\delta) &= \frac{\delta p r}{1-p} \\
Var(X_\delta)&= \frac{\delta^2 p r}{(1-p)^2}
\end{align*}
By requiring that these match the mean and variance of $X$, we find that:
\begin{align*}
r &= \frac{k\theta}{\theta-\delta} \\
q &= \frac{\theta-\delta}{\theta} 
\end{align*}
With these choices of $r$ and $q$, in the limit $\delta \rightarrow 0$ the discrete distribution converges to $\Gamma(k,\theta)$ in the following sense:
\begin{equation*}
\lim_{\delta \rightarrow 0} \frac{1}{\delta} b \left[\lfloor x /\delta \rfloor,r(k,\theta),q(k,\theta)\right] = \frac{x^{k-1}e^{-\frac{x}{\theta}}}{\Gamma(k) \theta^k}.
\end{equation*}
The advantage of thinking of the Gamma distribution as the limit of a negative binomial lies in the fact that the cascade distributions may be calculated explicitly in the discrete setting. We may extend the idea of non-integral convolution to the negative binomial distribution by making use of the following property. If $Y \sim NB(r,q)$ is a negative binomial variable, the the sum of $m$ independent copies of $Y$ has distribution $NB(mr,p)$. Setting $\delta=\tfrac{1}{m}$ we may think of our negative binomial approximation to the gamma distribution as the sum of $m$ negative binomial variables, $Y_i \in \{0,\delta, 2\delta, \ldots \}$, each with distribution:
\begin{equation*}
Y_i \sim NB \left( \frac{r}{m},q\right).
\end{equation*}
We will refer to this as the atomic distribution $A$. We can approximate non integral convolutions of the Gamma distribution as integral convolutions of the atomic distribution as follows: $\Gamma^{\ast X_n} \approx A^{\ast \lfloor m X_n \rfloor}$.

\subsection{Cascade distribution for the discrete state branching process}

From here on we set $k=2$ and $\theta=p$. The particular values of $r$ and $q$ in the atomic distribution required so that the negative binomial approximates the gamma distribution are:
\begin{align*}
r^* &= \frac{2 \delta p}{p-\delta} \\
q^* &= \frac{p-\delta}{p}
\end{align*}
Let $Z_\delta(m)$ be the total cascade size starting from $m$ individuals where the number of offspring produced by each individual is distributed according to the atomic distribution. We note that if $\delta = 1/m$ then $Z_\delta(m) \approx Z$. Because the the numbers of offspring produced by each individual are independent then $Z_\delta(m)$ has the same distribution as the sum of $m$ independent copies of $Z_\delta(1)$. Let $Y$ be an $ A$-distributed variate then:
\begin{equation}
\label{ZRec}
Z_\delta(1) = 1 + Z_\delta(Y)
\end{equation}
If $H(s)$ and $F(s)$ are the probability generating functions for $Z_\delta(1)$ and $Y$, then from equation (\ref{ZRec}) we have
\begin{align*}
H(s) &= \mathbf{E}(s^{1 + Z_\delta(Y)}) \\
&= s \mathbf{E}[ \mathbf{E}(s^{Z_\delta(Y)}\mid Y) ] \\
&= s \mathbf{E}[ (F(s))^Y ] \\
&= s F(H(s)).
\end{align*}
From the negative binomial mass function we have that:
\begin{equation*}
F(s) = \sum_{n=0}^\infty s^n b(n,r^*,q^*) = \left(\frac{1-q^*}{1-q^* s}\right)^{r^*}.
\end{equation*}
We are interested in the probability generating function of $Z_\delta(m)$, which is just $H^m(s)$. The coefficient of $s^n$ in this function may be determined using the Lagrange inversion formula:
\begin{align*}
[s^n]\left\{H^m(s)\right\} &= \frac{1}{n}[H^{n-1}]\left\{ \left(\frac{d}{dH} H^{m}\right) F^n(H)\right\} \\
&= \frac{m}{n} [H^{n-m}] F^n(H) \\
&= \frac{m}{n} \frac{\Gamma(n(1+r^*)-m)}{\Gamma(nr^*) \Gamma(n-m+1)}(1-q^*)^{r^*n} (q^*)^{n-m} \\
&= \PP\{Z_\delta(m) = n\}
\end{align*}
We now have the probability mass function for the cascade size in the discrete branching process which approximates the continuum process that we are interested in.

\subsection{The continuum limit}

We obtain the continuum cascade density function, which we will call $g(x)$, by setting $n=x/\delta$ and $m=1/\delta$ and then taking the limit $\delta \rightarrow 0$:
\begin{align}
g(x) &= \lim_{\delta \rightarrow 0} \frac {1}{\delta} \PP\{Z_\delta(m) = n\} \\
&= \frac{(x-1)^{2x-1}e^{-(\frac{1}{p} + 2 \ln p) x + \frac{1}{p}}}{x \Gamma(2x)}
\label{exact}
\end{align}
The asymptotic properties of $g(x)$ may be determined by making use of Stirling's approximation: $\Gamma(z+1) \sim \sqrt{2 \pi z} \left(\tfrac{z}{e}\right)^{z}$. The result is:
\begin{equation*}
g(x) \sim \frac{e^{\frac{1}{p}-2+\ln 2}}{2 \sqrt{\pi}} \frac{e^{-\left(\frac{1-2p}{p}+2 \ln 2p\right)x}}{x^{\frac{3}{2}}} \text{ as } x \rightarrow \infty
\end{equation*}
From this we see that the distribution is asymptotically a pure power law when $p=\tfrac{1}{2}$.

\subsection{Moments in the subcritical case}

Provided $p<\tfrac{1}{2}$, the distribution $g(x)$ is normalised and its moments are defined. It is useful to have explicit expressions for the mean and variance of $g(x)$ in this case. We first compute the mean and variance of the distribution of $Z_\delta(1)$ using the generating function relationship: $H(s)=sF(H(s))$, which after differentiation reveals that:
\begin{align*}
H'(s) &= \frac{F(H(s))}{1-s F'(H(s))} \\
H''(s) &= \frac{2 H'(s)F'(H(s))+s (H'(s))^2F''(H(s))}{1-sF'(s)}
\end{align*}
Using the expression for $F(s)$, together with the fact that when $p<\tfrac{1}{2}$ $H(1)=F(1)=1$, we find that:
\begin{align*}
\EE(Z_\delta(1)) = \delta H'(1) &= \frac{\delta}{1-2p} \\
\EE(Z_\delta(1)^2)- \EE(Z_\delta(1))^2 = \delta^2 (H''(1)+H'(1)-H'(1)^2) &= \frac{2 \delta p^2}{(1-2p)^3} 
\end{align*}
Since, in the limit $\delta \rightarrow 0$, $Z$ has then same distribution as the sum of $\tfrac{1}{\delta}$ copies of $Z_\delta(1)$, then:
\begin{align*}
\EE(Z) &= \frac{1}{1-2p} \\
\EE(Z^2)-\EE(Z)^2 &= \frac{2 p^2}{(1-2p)^3} 
\end{align*}
It is worth noting that these expressions would have been difficult to obtain by direct integration over $g(x)$.

\section{Probability of the event $\{Z<\infty\}$}

Numerical integration of the exact distribution (\ref{exact}) reveals that it is not normalized when $p>\tfrac{1}{2}$. The total probability weight is equal to $\PP\{Z<\infty\}$ which is less than one in the supercritical case. We will now show that:
\begin{equation}
\label{norm}
\int_1^\infty g(x) dx = e^{\chi(p)}
\end{equation}
where:
\begin{equation*}
\chi(p) = \left(2 W_{-1}\left(-\frac{e^{-\frac{1}{2p}}}{2p}\right)+\frac{1}{p}\right) \mathbf{1}_{[\frac{1}{2},\infty]}(p)
\end{equation*}
and $W_{-1}$ is one of the two real branches of the Lambert W function, the other being $W_0$. We may derive this result by mapping the discrete branching process on to a random walk, and then constructing a Martingale to which the optional sampling theorem may be applied.

We must show first that the total cascade size in a discrete branching process has the same distribution as a first passage time for a random walk. We consider the process with offspring distribution $A$ (the atomic distribution) and suppose that $X_0$ is the size of the zeroth generation. We identify $X_0$ with the initial position of the walker. The size, $X_1$, of the next generation is then the sum of $X_0$ $A$-distributed random variables. Now, suppose that $V \sim A$ and define $Q:=V-1$. We write the distribution of $Q$ as $A^-$, and note that $X_1$ has the same distribution as the position of a walker after $X_0$, $A^-$-distributed steps. Note that since $A^-$ has support $\{-1,0,1,2,\ldots \}$ then if the walker does reach the origin then it will do so at the $X_0$-th step. The sizes, $\{X_2,X_3,X_4,\ldots \}$ of subsequent generations may be viewed as the positions of the same walker after, respectively, $\{X_1,X_2,X_3,\ldots \}$ steps. The cascade ends when a generation has zero size, which occurs when the walker reaches the origin, and the total cascade size $X_0+X_1+X_2+\ldots$ is the total number of steps taken by the walker. The cascade size therefore has the same distribution as the time of first intersection of the $A^-$  distributed walker with the origin.

It remains to compute the probability that walker will return to the origin. We do this using martingales. The aim is to construct a martingale from the $A^-$ walk of the form:
\begin{equation*}
M_n = \alpha^{S_n}
\end{equation*}
where $S_n$ is the position of the walker after $n$ steps, and then to apply the optional sampling theorem, which states that under certain conditions (which will be satisfied for us) $\EE(M_T)=M_0$ when $T$ is a stopping time measurable w.r.t. the information contained in the walk up to the $n$th step. The value of $\alpha$ which makes $M_n$ a martingale satisfies the equation: $\EE(\alpha^Q)=1$, which has the explicit form:
\begin{equation*}
\frac{1}{\alpha}\left(\frac{1-q^*}{1-q^* \alpha}\right)^{r^*} = 1,
\end{equation*}
or, written in terms of $p$ and $\delta$:
\begin{equation}
\label{mart}
\left(\frac{\delta }{p}\right)^{\frac{2 p \delta
   }{p-\delta }} \left(\alpha  \left(\frac{\delta
   }{p}-1\right)+1\right)^{-\frac{2 p \delta }{p-\delta
   }}=\alpha
\end{equation}
Apart from the trivial solution $\alpha_1=1$, when $p>\tfrac{1}{2}$ this equation has another solution $\alpha_2 <1$ which approaches $1^-$ as $\delta \rightarrow 0$. Writing $\alpha=1-\epsilon$ and defining $x=\epsilon/\delta$ we find that the solutions to (\ref{mart}) converge, as $\delta \rightarrow 0$, to the solutions of
\begin{equation*}
x = 2 \ln(1+px),
\end{equation*}
the appropriate one being:
\begin{equation*}
x(p) =-2 W_{-1}\left(-\frac{e^{-\frac{1}{2p}}}{2p}\right)-\frac{1}{p}
\end{equation*}
To approximate the Gamma branching process we must start the atomic branching process with $X_0=1/\delta$. Optional sampling tells us that if $T$ is the first step at which the walk reaches the origin or $\infty$ then
\begin{equation*}
\EE(\alpha^T)= \PP\{Z_\delta <\infty\} = \alpha^{\frac{1}{\delta}}
\end{equation*}
Taking the limit $\delta \rightarrow 0$ we find that:
\begin{align*}
\PP\{Z<\infty\} &= \lim_{\delta \rightarrow 0}\PP\{Z_\delta <\infty\} \\
&= \lim_{\delta \rightarrow 0} (1-\delta x(p))^{\frac{1}{\delta}} \\
&= e^{-x(p)},
\end{align*}
which reproduces our claim (\ref{norm}).

\section{Concluding comment}

Branching processes are useful in the modelling of cascading failures, amongst many other applications. In the simplest case, the generations of the cascade are integer random variables, but this is not always the most appropriate model. The exact expression derived here for the probability distribution of cascade sizes the branching process with Gamma distributed generations is therefore likely to be useful tool in applications.

\end{document}